\documentclass{article}

\newtheorem{theorem}{Theorem}[section]
\newtheorem{corollary}[theorem]{Corollary}
\newtheorem{conjecture}[theorem]{Conjecture}

\newtheorem{proposition}[theorem]{Proposition}

 \evensidemargin\oddsidemargin
\begin{document}

\title{On the Roots of Independence Polynomials of Almost All Very Well-Covered
Graphs}
\author{Vadim E. Levit and Eugen Mandrescu \\
%EndAName
Department of Computer Science\\
Holon Academic Institute of Technology\\
52 Golomb Str., P.O. Box 305\\
Holon 58102, ISRAEL\\
\{levitv, eugen\_m\}@hait.ac.il}
\date{}
\maketitle

\begin{abstract}
If $s_{k}$ denotes the number of stable sets of cardinality $k$ in graph $G$%
, and $\alpha (G)$ is the size of a maximum stable set, then $%
I(G;x)=\sum\limits_{k=0}^{\alpha (G)}s_{k}x^{k}$ is the \textit{independence
polynomial} of $G$ (Gutman and Harary, 1983). A graph $G$ is \textit{%
very well-covered} (Favaron, 1982) if it has no isolated vertices,
its order equals $2\alpha (G)$ and it is \textit{well-covered}
(i.e., all its maximal independent sets are of the same size, M.
D. Plummer, 1970). For instance, appending a single pendant edge
to each vertex of $G$ yields a very well-covered graph, which we
denote by $G^{*}$. Under certain conditions, any well-covered
graph equals $G^{*}$ for some $G$ (Finbow, Hartnell and
Nowakowski, 1993).

The root of the smallest modulus of the independence polynomial of
any graph is real (Brown, Dilcher, and Nowakowski, 2000). The
location of the roots of the independence polynomial in the
complex plane, and the multiplicity of the root of the smallest
modulus are investigated in a number of articles

In this paper we establish formulae connecting the coefficients of $I(G;x)$
and $I(G^{*};x)$, which allow us to show that the number of roots of $I(G;x)$
is equal to the number of roots of $I(G^{*};x)$ different from $-1$, which
appears as a root of multiplicity $\alpha (G^{*})-\alpha (G)$ for $%
I(G^{*};x) $. We also prove that the real roots of $I(G^{*};x)$ are in $%
[-1,-1/2\alpha (G^{*}))$, while for a general graph of order $n$
we show that its roots lie in $\left| z\right| >1/(2n-1)$.

Hoede and Li (1994) posed the problem of finding graphs that can
be uniquely defined by their clique polynomials
(\textit{clique-unique graphs}). Stevanovic (1997) proved that
threshold graphs are clique-unique. Here, we demonstrate that the
independence polynomial distinguishes well-covered spiders
$(K_{1,n}^{*},n\geq 1)$ among well-covered trees.\newline

\textbf{keywords:}\textit{\ stable set, independence polynomial, root,
well-covered graph, clique-unique graph.}
\end{abstract}

\section{Introduction}

Throughout this paper $G=(V,E)$ is a simple (i.e., a finite, undirected,
loopless and without multiple edges) graph with vertex set $V=V(G)$ and edge
set $E=E(G).$ The \textit{complement} of $G$ is denoted by $\overline{G}$.
If $X\subset V$, then $G[X]$ is the subgraph of $G$ spanned by $X$. By $G-W$
we mean the subgraph $G[V-W]$, if $W\subset V(G)$. We also denote by $G-F$
the partial subgraph of $G$ obtained by deleting the edges of $F$, for $%
F\subset E(G)$, and we write shortly $G-e$, whenever $F$ $=\{e\}$. The
\textit{neighborhood} of a vertex $v\in V$ is the set $N_{G}(v)=\{w:w\in V$
\ \textit{and} $vw\in E\}$, and $N_{G}[v]=N_{G}(v)\cup \{v\}$; if there is
no ambiguity on $G$, we use $N(v)$ and $N[v]$, respectively. A vertex $v$ is
\textit{pendant} if its neighborhood contains only one vertex; an edge $e=uv$
is \textit{pendant} if one of its endpoints is a pendant vertex.

$K_{n},P_{n},C_{n},K_{n_{1},n_{2},...,n_{p}}$ denote respectively, the
complete graph on $n\geq 1$ vertices, the chordless path on $n\geq 1$
vertices, the chordless cycle on $n\geq 3$ vertices, and the complete
multipartite graph on $n_{1}+n_{2}+...+n_{p}$ vertices. As usual, a \textit{%
tree} is an acyclic connected graph. A \textit{spider} is a tree having at
most one vertex of degree $\geq 3$, \cite{HedLaskar}.

A \textit{stable} set in $G$ is a set of pairwise non-adjacent vertices. A
stable set of maximum size will be referred to as a \textit{maximum stable
set} of $G$, and the \textit{stability number }of $G$, denoted by $\alpha
(G) $, is the cardinality of a maximum stable set in $G$. Let $s_{k}$ be the
number of stable sets in $G$ of cardinality $k$. The polynomial
\[
I(G;x)=\sum\limits_{k=0}^{\alpha (G)}s_{k}x^{k}
\]
is called the \textit{independence polynomial} of $G$, (Gutman and Harary,
\cite{GutHar}), or the \textit{clique polynomial} of the complement of $G$
(Hoede and Li, \cite{HoedeLi}).

While further we will follow the notation of Gutman and Harary, it is worth
mentioning that in \cite{FisherSolow} the \textit{dependence} \textit{%
polynomial} $D(G;x)$ of a graph $G$ is defined as $D(G;x)=I(\overline{G}%
;-x)=\sum\limits_{k=0}^{\omega \left( G\right) }(-1)^{k}s_{k}x^{k},\omega
\left( G\right) =\alpha (\overline{G})$ , where $s_{k}$ is the number of
stable sets of size $k$ in $\overline{G}$. In \cite{GoldSantini}, $D(G;x)$
is defined as the \textit{clique polynomial} of $G$. In \cite{DohPonTit},
the independence polynomial appears as a particular case of a two-variable
generalized chromatic polynomial.

A graph $G$ is called \textit{well-covered} if all its maximal stable sets
are of the same cardinality, (Plummer, \cite{Plum}). If, in addition, $G$
has no isolated vertices and its order equals $2\alpha (G)$, then $G$ is
\textit{very well-covered} (Favaron, \cite{Fav1}).

Throughout this paper, by $G^{*}$ we mean the graph obtained from $G$ by
appending a single pendant edge to each vertex of $G$, \cite
{DuttonChanBrigham}. In \cite{ToppLutz}, $G^{*}$ is denoted by $G\circ K_{1}$
and is defined as the \textit{corona} of $G$ and $K_{1}$. We refer to $G$ as
to a \textit{skeleton} of $G^{*}$. Let us remark that $G^{*}$ is
well-covered (see, for instance, \cite{LevMan0}), and $\alpha (G^{*})=n$. In
fact, $G^{*}$ is very well-covered. Moreover, the following result shows
that, under certain conditions, any well-covered graph equals $G^{*}$ for
some graph $G$.

\begin{theorem}
\cite{FinHarNow}\label{th3} Let $G$ be a connected graph of girth $\geq 6$,
which is isomorphic to neither $C_{7}$ nor $K_{1}$. Then $G$ is well-covered
if and only if its pendant edges form a perfect matching.
\end{theorem}

In other words, Theorem \ref{th3} shows that apart from $K_{1}$ and $C_{7}$,
connected well-covered graphs of girth $\geq 6$ are very well-covered. In
particular, a tree $T$ is well-covered if and only if $T=K_{1}$ or it has a
perfect matching consisting of pendant edges (Ravindra, \cite{Ravindra}). It
turns out that a tree $T\neq K_{1}$ is well-covered if and only if it is
very well-covered. An alternative characterization of well-covered trees is
the following:

\begin{theorem}
\cite{LevitMan1} A tree $T$ is well-covered if and only if either $T$ is a
well-covered spider, or $T$ is obtained from a well-covered tree $T_{1}$ and
a well-covered spider $T_{2}$, by adding an edge joining two non-pendant
vertices of $T_{1},T_{2}$,\ respectively.
\end{theorem}

\begin{figure}[h]
\setlength{\unitlength}{0.8cm}
\begin{picture}(5,3)\thicklines

  \put(2,0){\circle*{0.29}}
  \put(2,0.5){\makebox(0,0){$K_{1}$}}
  \put(3.5,0){\circle*{0.29}}
  \put(3.5,1){\circle*{0.29}}
  \put(3.5,0){\line(0,1){1}}
  \put(3,0.5){\makebox(0,0){$K_{2}$}}

  \multiput(5,0)(1,0){2}{\circle*{0.29}}
  \multiput(5,1)(1,0){2}{\circle*{0.29}}
  \put(5,0){\line(1,0){1}}
  \multiput(5,0)(1,0){2}{\line(0,1){1}}
  \put(4.5,0.5){\makebox(0,0){$P_{4}$}}

  \multiput(8,0.5)(3,0){3}{\circle*{0.29}}
  \multiput(8,1.5)(1,0){7}{\circle*{0.29}}
  \multiput(9,2.5)(1,0){2}{\circle*{0.29}}
  \multiput(12,2.5)(1,0){2}{\circle*{0.29}}
  \multiput(8,0.5)(1,0){6}{\line(1,0){1}}
  \multiput(8,0.5)(3,0){3}{\line(0,1){1}}
  \multiput(9,1.5)(1,0){2}{\line(0,1){1}}
  \multiput(12,1.5)(1,0){2}{\line(0,1){1}}
  \put(11,0.5){\line(-2,1){2}}
  \put(11,0.5){\line(-1,1){1}}
  \put(11,0.5){\line(1,1){1}}
  \put(11,0.5){\line(2,1){2}}
  \put(11,0){\makebox(0,0){$b_{0}$}}
 \put(11,1.9){\makebox(0,0){$a_{0}$}}
 \put(8,0){\makebox(0,0){$b_{1}$}}
 \put(7.6,1.7){\makebox(0,0){$a_{1}$}}
 \put(9,1){\makebox(0,0){$b_{2}$}}
 \put(8.6,2.5){\makebox(0,0){$a_{2}$}}
 \put(10.45,1.6){\makebox(0,0){$b_{3}$}}
 \put(10.5,2.5){\makebox(0,0){$a_{3}$}}
 \put(12.45,1.6){\makebox(0,0){$b_{4}$}}
 \put(12.5,2.5){\makebox(0,0){$a_{4}$}}
 \put(13.1,1){\makebox(0,0){$b_{5}$}}
 \put(13.5,2.5){\makebox(0,0){$a_{5}$}}
 \put(14,0){\makebox(0,0){$b_{6}$}}
 \put(14.5,1.5){\makebox(0,0){$a_{6}$}}
 \put(7.2,0.9){\makebox(0,0){$S_{6}$}}

 \end{picture}
\caption{Well-covered spiders.}
\label{fig99}
\end{figure}
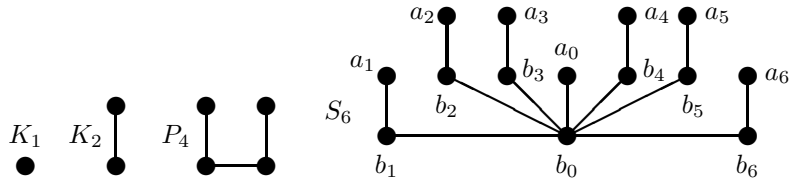

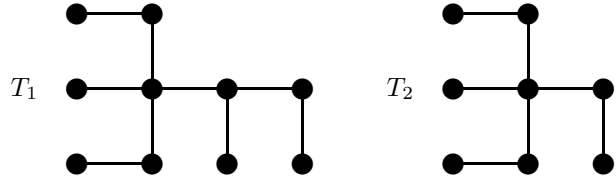
\begin{figure}[h]
\setlength{\unitlength}{1cm}%
\begin{picture}(5,2.5)\thicklines

  \multiput(3,0)(1,0){4}{\circle*{0.29}}
  \multiput(3,1)(1,0){4}{\circle*{0.29}}
  \multiput(3,2)(1,0){2}{\circle*{0.29}}
  \put(3,0){\line(1,0){1}}
  \put(3,1){\line(1,0){3}}
  \put(3,2){\line(1,0){1}}
  \put(4,0){\line(0,1){2}}
  \put(5,0){\line(0,1){1}}
  \put(6,0){\line(0,1){1}}

\put(2.3,1){\makebox(0,0){$T_{1}$}}

  \multiput(8,0)(1,0){3}{\circle*{0.29}}
  \multiput(8,1)(1,0){3}{\circle*{0.29}}
  \multiput(8,2)(1,0){2}{\circle*{0.29}}
  \put(8,0){\line(1,0){1}}
  \put(8,1){\line(1,0){2}}
  \put(8,2){\line(1,0){1}}
  \put(9,0){\line(0,1){2}}
  \put(10,0){\line(0,1){1}}

  \put(7.3,1){\makebox(0,0){$T_{2}$}}

 \end{picture}
\caption{Two (very) well-covered trees.}
\label{fig5}
\end{figure}

The roots of independence polynomials of (well-covered) graphs are not
necessarily real, even if they are trees. For instance, the trees $%
T_{1},T_{2}$ in Figure \ref{fig5} are very well-covered, their independence
polynomials are respectively,
\begin{eqnarray*}
I(T_{1};x)
&=&(1+x)^{2}(1+2x)(1+6x+7x^{2})=1+10x+36x^{2}+60x^{3}+47x^{4}+14x^{5}, \\
I(T_{2};x) &=&(1+x)(1+7x+14x^{2}+9x^{3})=1+8x+21x^{2}+23x^{3}+9x^{4},
\end{eqnarray*}
but only $I(T_{1};x)$ has all the roots real.
Moreover, it is easy to check that the complete $n$-partite graph $%
G=K_{\alpha ,\alpha ,...,\alpha }$ is well-covered, $\alpha (G)=\alpha $,
and its independence polynomial $I(G;x)=n(1+x)^{\alpha }-(n-1)$ has only one
real root, whenever $\alpha $ is odd, and exactly two real roots, for any
even $\alpha \geq 2$.

The roots of the independence polynomial of (well-covered) graphs are
investigated in a number of papers, as \cite{BrownDilNow}, \cite{BrownNow},
\cite{FisherSolow}, \cite{GoldSantini}, \cite{MehHaj}. Denoting by $\xi
_{\min },\xi _{\max }$ the smallest and the largest real root of $I(G;x)$,
respectively, we get that $\xi _{\min }\leq \xi _{\max }<0$, since all the
coefficients of $I(G;x)$ are positive. Let us recall the following known
results.

\begin{proposition}
\label{prop1}If $G$ is a graph of order $n\geq 2$, then:

\emph{(i)} \cite{FisherSolow} the smallest (in absolute value) root $\lambda
$ of $I(\overline{G};-x)$

\qquad satisfies $0<\lambda \leq \alpha (G)/n$, i.e., $-\frac{\alpha \left(
G\right) }{n}\leq \xi _{\max }<0$;

\emph{(ii)} \cite{GoldSantini} $I(\overline{G};-x)$ has only one root of
smallest modulus $\rho $ and, furthermore,

$\qquad 0<\rho \leq 1$, i.e., $\xi _{\max }$ is unique and $0<\left| \xi
_{\max }\right| \leq 1$;

\emph{(iii)} \cite{BrownDilNow} a root of smallest modulus of $I(G;x)$ is
real, for any graph $G$,

\qquad i.e., for $I(G;x)$ there exists $\xi _{\max }$;

\emph{(iv)} \cite{BrownDilNow} for a well-covered graph $G$ on $n\geq 1$
vertices, the roots of $I(G;x)$ lie in

\qquad the annulus $1/n\leq \left| z\right| \leq \alpha (G)$,

\qquad furthermore, there is a root on the boundary if and only if $G$ is
complete;

\emph{(v)} \cite{MehHaj} if $\mu $ is the greatest real root of $I(\overline{%
G};x)$, then $\alpha (\overline{G})\leq -1/\mu $,

\qquad i.e., $-1/\alpha (\overline{G})\leq \xi _{\max }$.
\end{proposition}

It is also shown in \cite{BrownDilNow} that for any well-covered graph $G$
there is a well-covered graph $H$ with $\alpha (G)=\alpha (H)$ such that $G$
is an induced subgraph of $H$ and $I(H;x)$ has all its roots simple and
real. In \cite{BrownNow} the problem of determining the maximum modulus of
roots of independence polynomials for fixed stability number is completely
solved, namely, the bound is $(n/\alpha )^{\alpha -1}+O(n^{\alpha -2})$,
where $\alpha =\alpha (G)$ and $n=\left| V(G)\right| $.

Let us mention that there are non-isomorphic (well-covered) graphs with the
same independence polynomial (see, for example, Figures \ref{fig63}, \ref
{fig70}). Following Hoede and Li, \cite{HoedeLi}, $G$ is called a \textit{%
clique-unique graph} if the relation $I(\overline{G};x)=I(\overline{H};x)$
implies that $\overline{G}$ and $\overline{H}$ are isomorphic (or,
equivalently, $G$ and $H$ are isomorphic). One of the problems they proposed
was to determine clique-unique graphs (Problem 4.1, \cite{HoedeLi}). In \cite
{Stevan}, Stevanovic proved that the \textit{threshold graphs} (i.e., graphs
having no induced subgraph isomorphic to a either a $P_{4}$, or a $C_{4}$,
or a $\overline{C_{4}}$, defined by Chvatal and Hammer, \cite{ChvatHam}) are
clique-unique graphs.

In this paper we emphasize a number of formulae transforming the
coefficients of $I(G;x)$ to the coefficients of $I(G^{*};x)$, and vice
versa. Based on these results, we deduce some properties connecting $I(G;x)$
and $I(G^{*};x)$. For instance, it is shown that the number of roots of $%
I(G;x)$ is equal to the number of roots of $I(G^{*};x)$ different from $-1$.
Moreover, $-1$ is a $\left( \alpha (G^{*})-\alpha (G)\right) $-folded root
of $I(G^{*};x)$.

We also strengthen\emph{\ }Proposition \ref{prop1}\emph{(iv)}, as it
concerns the real roots. Namely, we prove that the real roots of the
independence polynomial of a non-complete well-covered graph $G,G\neq C_{7}$
and of girth $\geq 6$, are in $[-1,-1/n)$, where $n=2\alpha (G)$.

As an application of our findings, we show that independence polynomials
distinguish between well-covered spiders and general well-covered trees.

\section{The polynomials $I(G;x),I(G^{*};x)$ and their roots}

As we saw in the introduction, the skeleton $G=(V,E),V=\{v_{i}:1\leq i\leq
n\}$ defines $G^{*}$ using a set of additional vertices $U=\{u_{i}:1\leq
i\leq n\}$ as follows
\[
G^{*}=(V\cup U,E\cup \{u_{i}v_{i}:1\leq i\leq n\}).
\]

Let us denote the independence polynomials of $G$ and $G^{*}$ as
\[
I(G;x)=\sum\limits_{k=0}^{\alpha (G)}s_{k}x^{k} \ and \
I(G^{*};x)=\sum\limits_{k=0}^{\alpha (G^{*})}t_{k}x^{k},
\]
respectively.

\begin{theorem}
\label{th2}For any graph $G$ of order $n$ the following assertions are true:

\emph{(i)} the independence polynomial of $G^{*}$ is
\[
I(G^{*};x)=\sum\limits_{k=0}^{\alpha (G)}s_{k}\cdot x^{k}\cdot
(1+x)^{n-k}=(1+x)^{\alpha (G^{*})-\alpha (G)}\cdot \sum\limits_{k=0}^{\alpha
(G)}s_{k}\cdot x^{k}\cdot (1+x)^{\alpha (G)-k};
\]

\qquad and the formulae connecting the coefficients of $I(G;x)$ and $%
I(G^{*};x)$ are
\begin{eqnarray*}
t_{k} &=&\sum\limits_{j=0}^{k}s_{j}\cdot {n-j \choose n-k}, k\in
\{0,1,...,\alpha (G^{*})=n\}, \\
s_{k} &=&\sum\limits_{j=0}^{k}(-1)^{k+j}\cdot t_{j}\cdot {n-j \choose n-k}%
,k\in \{0,1,...,\alpha (G)\},
\end{eqnarray*}

\qquad for example, $t_{0}=1$ and $t_{n}=s_{0}+s_{1}+...+s_{\alpha (G)}$
(the fact that the number of

\qquad stable sets of $G$ equals the highest coefficient of $I(G^{*};x)$ is
mentioned in an

\qquad implicit form in \cite{DuttonChanBrigham});

\emph{(ii)} $t_{0}\leq t_{1}\leq ...\leq t_{j}$, where $j=\left\lceil
n/2\right\rceil $.
\end{theorem}

\setlength {\parindent}{0.0cm}\textbf{Proof.} It is easy to observe that $%
\alpha (G^{*})=\left| V\right| =n$, and, correspondingly, $%
I(G^{*};x)=\sum\limits_{k=0}^{n}t_{k}x^{k}$.\setlength
{\parindent}{3.45ex}

\emph{(i)} Clearly, $t_{0}=s_{0}=1$. A stable set $S$ in $H$ of size $%
m,1\leq m\leq n$, can be obtained as follows:

\begin{itemize}
\item  $S\subseteq V$ , only for $m\leq \alpha (G)$, and there are $s_{m}$
sets of this kind, or

\item  $S\subseteq U$, and the number of stable sets of this form
is ${n \choose n-m}$, or

\item  $S=S_{1}\cup S_{2}$ with $S_{1}\subseteq V,\left|
S_{1}\right| =j\leq \alpha (G),S_{2}\subseteq U-\{u_{i}:v_{i}\in
S_{1}\},\left| S_{2}\right| =m-j $, and there exist $s_{j}$ sets
of the form $S_{1}$ and ${n-j \choose n-m}$ sets of the form
$S_{2}$ (because $\left| U-\{u_{i}:v_{i}\in S_{1}\}\right|
=n-j$); therefore, there are ${n-j \choose n-m}\cdot s_{j}$ stable sets in $%
G^{*}$ of this kind.
\end{itemize}

Consequently, we infer that
\[
t_{m}=\sum\limits_{j=0}^{m}{n-j \choose n-m}\cdot s_{j}=\sum\limits_{j=0}^{m}%
{n-j \choose m-j}\cdot s_{j},
\]
where, clearly, $s_{j}=0\ $for$\ j>\alpha (G)$. On the other hand, it is
easy to see that the coefficient of $x^{m},1\leq m\leq n$, in the polynomial
\[
\sum\limits_{k=0}^{\alpha (G)}s_{k}\cdot x^{k}\cdot
(1+x)^{n-k}=\sum\limits_{k=0}^{\alpha
(G)}(\sum\limits_{j=0}^{n-k}s_{k}\cdot {n-k \choose j}\cdot
x^{k+j})
\]
is exactly $t_{m}$. Therefore, the equality $I(G^{*};x)=\sum\limits_{k=0}^{%
\alpha (G)}s_{k}\cdot x^{k}\cdot (1+x)^{n-k}$ is true.

A proof for the inverse formulae
\[
s_{k}=\sum\limits_{j=0}^{k}(-1)^{k+j}\cdot t_{j}\cdot {n-j \choose
n-k},k\in \{0,1,...,\alpha (G)\},
\]
can be found in \cite{Riordan} and \cite{StantonSprott}.

\emph{(ii)} For $n\geq 3$, let us observe that ${n \choose 1}\leq
{n \choose 2} $ and ${n-1 \choose 0}\leq {n-1 \choose 1}$ imply
\[
t_{0}=s_{0}\leq {n \choose 1}\cdot s_{0}+{n-1 \choose 0}\cdot
s_{1}=t_{1}\leq {n \choose 2}\cdot s_{0}+{n-1 \choose 1}\cdot
s_{1}+{n-2 \choose 0}\cdot s_{2}=t_{2}.
\]
Since, in general, ${n \choose 0}\leq {n \choose 1}\leq ...\leq {n \choose %
\left\lceil n/2\right\rceil }$ is true for the binomial coefficients, we
deduce that for $i+1\leq \left\lceil n/2\right\rceil $ we have:
\begin{eqnarray*}
t_{i} &=&{n \choose i}\cdot s_{0}+{n-1 \choose i-1}\cdot s_{1}+...+{n-i \choose %
0}\cdot s_{i} \\
&\leq &{n \choose i+1}\cdot s_{0}+{n-1 \choose i}\cdot s_{1}+...+{n-i \choose 1}%
\cdot s_{i}+{n-i-1 \choose 0}\cdot s_{i+1}=t_{i+1}.
\end{eqnarray*}
Therefore, we may conclude that $t_{0}\leq t_{1}\leq ...\leq t_{j}$, where $%
j=\left\lceil n/2\right\rceil $. \rule{2mm}{2mm}\newline

Actually, the inequality from Theorem \ref{th2}\emph{(ii)} is true for any
well-covered graph. We infer this fact (Proposition \ref{prop2}\emph{(iii)})
as a simple consequence of a generalization of the well-known Theorem of
Euler (Proposition \ref{prop2}\emph{(i)}), stating that
\[
\sum\limits_{v\in V(G)}\deg (v)=2\left| E(G)\right| ={2 \choose
1}\left| E(G)\right| .
\]
Let $Q_{i}$ be an $i$\textit{-clique} in a graph $G$, i.e., a clique of size
$i$ in $G$; by $\deg _{j}(Q_{i})$ we mean the number of cliques of size $%
j\geq i$ that contains $Q_{i}$. In particular, for an $1$-clique, say $\{v\}$%
, $\deg _{2}(\{v\})$ equals the usual degree of the vertex $v$.

\begin{proposition}
\label{prop2}\emph{(i)} The equality
\[
\sum \{\deg _{j}(Q_{i}):Q_{i}\ is\ an\ i-clique\ in\
\overline{G}\}={j \choose i}\cdot s_{j}
\]
is true for any graph $G$.

\emph{(ii)} If $G$ is a well-covered graph and $1\leq i\leq j\leq
\alpha =\alpha (G)$, then ${\alpha -i \choose j-i}\cdot s_{i}\leq
{j \choose i}\cdot s_{j}$.

\emph{(iii)} If $G$ is a well-covered graph, then $s_{k-1}\leq s_{k}$ for
any $1\leq k\leq (\alpha (G)-1)/2$.
\end{proposition}

\setlength {\parindent}{0.0cm}\textbf{Proof.} \emph{(i)} Any $j$-clique
includes ${j \choose i}$ cliques of size $i\leq j$, and the number of $j$%
-cliques in $\overline{G}$ is exactly $s_{j}$. Consequently, there are $%
{j \choose i}\cdot s_{j}$ different inclusions $Q_{i}\subseteq
Q_{j}$, where $Q_{i}$ and $Q_{j}$ are an $i$-clique and a
$j$-clique, correspondingly. Since, according to the definition,
$\deg _{j}(Q_{i})$ is equal to the number of
cliques of size $j\geq i$ containing $Q_{i}$, the proof is complete.%
\setlength
{\parindent}{3.45ex}

\emph{(ii)} Since $G$ is well-covered, any $i$-clique $Q_{i}$ of $\overline{G%
}$ is included in an $\alpha $-clique $Q_{\alpha }$ of $\overline{G}$, and
there are ${\alpha -i \choose j-i}$ cliques of size $j$ in the clique $%
Q_{\alpha }$ that contains $Q_{i}$. Hence,
\[
{\alpha -i \choose j-i}\leq \deg _{j}(Q_{i}).
\]

Taking into account that the number of $i$-cliques of
$\overline{G}$ is exactly $s_{i}$, and using already proved
Proposition \ref{prop2}\emph{(i)}, we obtain

\[
{\alpha -i \choose j-i}\cdot s_{i}\leq \sum \{\deg _{j}(Q_{i}):Q_{i}\ is\ %
an\ i-clique\ in \overline{G}\}={j \choose i}\cdot s_{j},
\]
which completes the proof.

\emph{(iii)} Substituting $j=k$ and $i=k-1$ in Proposition \ref{prop2}\emph{%
(ii)}, we infer that
\[
(\alpha -k+1)\cdot s_{k-1}\leq k\cdot s_{k},
\]
which further leads to $s_{k-1}\leq s_{k}$, whenever $\alpha -k+1\geq k$,
i.e., $k\leq (\alpha +1)/2$. \rule{2mm}{2mm}\newline

Let us remark that Proposition \ref{prop2}\emph{(ii) }strengthens one
assertion from \cite{BrownDilNow}, where for any well-covered graph $G$ on $%
n $ vertices it is proved that $s_{k-1}\leq k\cdot s_{k}$ and also
$s_{k}\leq (n-k+1)\cdot s_{k-1},1\leq k\leq \alpha (G)$.

After finishing this paper we found that the statements contained in
Proposition \ref{prop2}$\emph{(ii),(iii)}$ were shown independently, in \cite
{MichaelTraves}.

\begin{theorem}
\label{th4}For any graph $G$ of order $n$ and with at least one
edge, the following assertions are true:

\emph{(i)} $G^{*}$ has an even number of stable sets;

\qquad moreover, the number of stable sets of $G^{*}$ is divisible by $%
2^{n-\alpha (G)}$;

\emph{(ii)} if $x\notin \{-1,0\}$, then $x^{n}\cdot
I(G^{*};1/x)=(1+x)^{n}\cdot I(G;1/(1+x))$;

\qquad substituting $x$ by $1/x$ one gets $I(G^{*};x)=(1+x)^{n}\cdot
I(G;x/(1+x))$;

\qquad further, $I(G^{*};x-1)=x^{n}\cdot I(G;1-1/x))$ is obtained by
changing $x$ into $x-1$;

\qquad for instance, $x=-1$ gives $I(G^{*};-2)=(-1)^{n}\cdot I(G;2)$.

\emph{(iii)} there exists a bijection between the set of roots of $I(G^{*};x)
$ different from $-1$

\qquad and the set of roots of $I(G;x)$, respecting the multiplicities of
the roots;

\qquad moreover, rational roots correspond to rational roots, and real roots

\qquad correspond to real roots;

\emph{(iv) }$-1$ is a root of $I(G^{*};x)$ with the multiplicity $\alpha
(G^{*})-\alpha (G)\geq 1$;

\emph{(v)} for any positive integer $k$ there exists a well-covered tree $%
H_{k}$,

\qquad such that $I(H_{k};-1/k)=0$;

\emph{(vi)} if $x<-1$, then $I(G^{*};x)\neq 0$, moreover, if $n$ is odd,
then $I(G^{*};x)<0$,

\qquad while for $n$ even, $I(G^{*};x)>0$.
\end{theorem}

\setlength {\parindent}{0.0cm}\textbf{Proof.} As in Theorem \ref{th2}, $%
\alpha (G^{*})=\left| V\right| =n$, and $I(G^{*};x)=\sum%
\limits_{k=0}^{n}t_{k}x^{k}$.\setlength
{\parindent}{3.45ex}

\emph{(i)} By Theorem \ref{th2}\emph{(i)}, it follows that $%
I(G^{*};1)=2^{n-\alpha (G)}\cdot \sum\limits_{k=0}^{\alpha
(G)}s_{k}\cdot 2^{\alpha (G)-k}$. Hence,
$I(G^{*};1)=t_{0}+t_{1}+...+t_{n}$ is a positive integer divisible
by $2^{n-\alpha (G)}=2^{\alpha (G^{*})-\alpha (G)}\geq 2$, because
$E(G) \neq \emptyset $ ensures that $\alpha (G)<n$.

\emph{(ii)} The equality $I(G^{*};x)=\sum\limits_{k=0}^{n}s_{k}\cdot
x^{k}\cdot (1+x)^{n-k}$ from Theorem \ref{th2}\emph{(i) }implies that
\begin{eqnarray*}
I(G^{*};1/x) &=&\sum\limits_{k=0}^{n}s_{k}\cdot 1/x^{k}\cdot (1+1/x)^{n-k}=
\\
&=&(1+x)^{n}/x^{n}\cdot \sum\limits_{k=0}^{n}s_{k}\cdot
[1/(1+x)]^{k}=(1+x)^{n}/x^{n}\cdot I(G;1/(1+x)),
\end{eqnarray*}
which can be written as $x^{n}\cdot I(G^{*};1/x)=(1+x)^{n}\cdot I(G;1/(1+x))$%
.

\emph{(iii)} Let $A=\{x:I(G;x)=0\}$ and $B=\{x:I(G^{*};x)=0,x\neq -1\}$.

Changing $x$ into $1/x-1$ in $x^{n}\cdot I(G^{*};1/x)=(1+x)^{n}\cdot
I(G;1/(1+x))$, we obtain $I(G;x)=(1-x)^{n}\cdot I(G^{*};x/(1-x))$ which
shows that there is an injection
\[
f_{1}:A\rightarrow B,\ f_{1}(x)=\frac{x}{1-x}
\]
from the set of roots of $I(G;x)$ to the set of the roots of $I(G^{*};x)$
different from $-1$, because $1$ can not be a root of $I(G;x)$ and $%
x/(1-x)\neq -1$.

Using $I(G^{*};x)=(1+x)^{n}\cdot I(G;x/(1+x))$, we see that there is an
injection
\[
f_{2}:B\rightarrow A,\ f_{2}(x)=\frac{x}{1+x}
\]
from the set of the roots of $I(G^{*};x)$ different from $-1$ to the set of
the roots of $I(G;x)$.

Together these claims give us a bijection $f=f_{1}=f_{2}^{-1}$
between the sets $A$ and $B$. Clearly, this bijection respects
belonging of roots to any subfield of $\mathbf C $, for instance,
for $\mathbf Q $, $\mathbf R $, etc.

Further, we get
\[
I^{\prime }(G^{*};x)=(1+x)^{n-2}\cdot [n\cdot (1+x)\cdot
I(G;x/(1+x))+I^{\prime }(G;x/(1+x))],
\]
which assures that if $b\in B$ has the multiplicity $m(b)=2$, i.e., $%
I(G^{*};b)=I^{\prime }(G^{*};b)=0$, then for $a=f^{-1}(b)$ we obtain $%
I(G;a)=I^{\prime }(G;a)=0$, that is $a=b/(1+b)$ must be a root of $I(G;x)$
of multiplicity $m(a)=2$, at least. Similarly,
\begin{eqnarray*}
I^{\prime \prime }(G^{*};x) &=&(1+x)^{n-3}\cdot [n\cdot (n-1)\cdot
(1+x)\cdot I(G;x/(1+x))+ \\
&&+2(n-1)\cdot I^{\prime }(G;x/(1+x))+I^{\prime }(G;x/(1+x))/(1+x)],
\end{eqnarray*}
ensures that if $b\in B$ has multiplicity $m(b)=3$, i.e., $%
I(G^{*};b)=I^{\prime }(G^{*};b)=I^{\prime \prime }(G^{*};b)=0$, then $%
a=b/(1+b)=f^{-1}(b)$ must be a root of $I(G;x)$ of multiplicity $m(a)=3$, at
least, because $I^{\prime \prime }(G^{*};b)=0=I(G^{*};b)=I^{\prime
}(G^{*};b) $ implies also $I^{\prime \prime }(G;a)=0$.

In this way we deduce that any root $b\in B$ leads to a root $a=f^{-1}(b)\in
A$ of multiplicity $m(a)\geq m(b)$.

Similarly, using the relation $I(G;x)=(1-x)^{n}\cdot I(G^{*};x/(1-x))$ we
infer that any root $a\in A$ gives rise to a root $b=f(a)\in B$ of
multiplicity $m(b)\geq m(a)$.

Thus, $m(f\left( a\right) )=m(a)$, for any $a\in A$. In other words, the
bijection $f$ respects the multiplicities of the roots.

\emph{(iv) }The equality
\[
I(G^{*};x)=\sum\limits_{k=0}^{n}s_{k}\cdot x^{k}\cdot
(1+x)^{n-k}=(1+x)^{\alpha (G^{*})-\alpha (G)}\cdot \sum\limits_{k=0}^{\alpha
(G)}s_{k}\cdot x^{k}\cdot (1+x)^{\alpha (G)-k},
\]
implies that $-1$ is a root of $I(G^{*};x)$ with the multiplicity at least $%
n-\alpha (G)=\alpha (G^{*})-\alpha (G)\geq 1$, where the
inequality goes from the hypothesis that $G$\ has at least one
edge.

On the other hand, using Theorem \ref{th4}\emph{(iii) }we get the following equality $%
\sum\limits_{a\in A}m(a)=\sum\limits_{b\in B}m(b)$. Since the polynomials $%
I(G;x)$ and $I(G^{*};x)$ are of degrees $\alpha (G)$ and $\alpha (G^{*})$,
respectively, the definitions of the sets $A$ and $B$ immediately give
\[
\alpha (G)=\sum\limits_{a\in A}m(a),\alpha (G^{*})=m(-1)+\sum\limits_{b\in
B}m(b),
\]
which finally provide the exact value of the multiplicity of $-1$ in $%
I(G^{*};x)$, namely, $m(-1)=\alpha (G^{*})-\alpha (G).$

\emph{(v)} Let $T\neq K_{1}$ be some tree, and $H_{1}=T^{*}$. $H_{1}$ is
well-covered and, according to Theorem \ref{th4}\emph{(iv)},\emph{\ }we get that $%
I(H_{1};-1)=0$. Let now $H_{1}$ be the skeleton of $H_{2}$. Taking
$x=1/2$ in the relation $I(H_{2};x-1)=x^{n}\cdot
I(H_{1};1-1/x)$,\emph{\ }we infer that $I(H_{2};-1/2)=1/2^{n}\cdot
I(H_{1};-1)=0$. If $H_{3}=H_{2}^{*}$, then
for $x=2/3$ in $I(H_{3};x-1)=x^{n}\cdot I(H_{2};1-1/x)$, we obtain $%
I(H_{3};-1/3)=(2/3)^{n}\cdot I(H_{2};-1/2)=0$. In general, if $H_{k-1}$ is
the skeleton of $H_{k}$, then taking $x=(k-1)/k$ in $I(H_{k};x-1)=x^{n}\cdot
I(H_{k-1};1-1/x)$, it implies
\[
I(H_{k};-1/k)=((k-1)/k)^{n}\cdot I(H_{k-1};-1/(k-1))=0
\]
and, clearly, $H_{k}$ is well-covered.

\emph{(vi)} The equality $I(G^{*};x)=(1+x)^{n}\cdot I(G;x/(1+x))$ from \emph{%
(ii)} shows also that $I(G^{*};x)\neq 0$ for any $x<-1$, since in
this case, $x/(1+x)>0$ and $I(G;x/(1+x))>0$, as well. Clearly, if
$n$ is odd, then $I(G^{*};x)<0$ for any $x<-1$, while for $n$
even, $I(G^{*};x)>0$ for any $x<-1$. \rule{2mm}{2mm}

\begin{corollary}
The number of stable sets of any well-covered tree $\neq K_{2}$ is divisible
by some power of $2$, while there are trees having an odd number of stable
sets; $K_{2}$ is the unique well-covered tree with an odd number of stable
sets.
\end{corollary}

\setlength {\parindent}{0.0cm}\textbf{Proof.} Let $T$ be a well-covered
tree. Clearly, $K_{1}$ has two stable sets, and $K_{2}$ has three stable
sets. If $T\neq K_{1},K_{2}$, then, according to Ravindra's result, $T$ has
a perfect matching consisting of pendant edges, i.e., $T=G^{*}$ for some
tree $G$. Then, according to Theorem \ref{th4}\emph{(i)}, $%
I(G^{*};1)=t_{0}+t_{1}+...+t_{n}=I(T;1)$ is a positive integer number
divisible by $2^{n-\alpha (G)}$. In other words, the number of stable sets
of $T$ is divisible by some power of $2$. However, $I(P_{3};x)=1+3x+x^{2}$
implies $I(P_{3};1)=5$, i.e., $P_{3}$ has an odd number of stable sets. On
the other hand, $I(G_{3};x)=1+6x+10x^{2}+6x^{3}+x^{4}$ gives $I(G_{3};1)=24$%
, (where $G_{3}$ is depicted in Figure \ref{fig63}), i.e., there are
non-well-covered trees having an even number of stable sets. \rule{2mm}{2mm}%
\setlength
{\parindent}{3.45ex}\newline

As a simple application of Theorem \ref{th2}\emph{(i)}, let us
notice that for any $n\geq 1,I(K_{n};x)=1+nx,\alpha
(K_{n}^{*})=n$, and therefore,
\[
I(K_{n}^{*};x)=(1+x)^{n-1}\cdot \sum\limits_{k=0}^{1}s_{k}\cdot x^{k}\cdot
(1+x)^{1-k}=(1+x)^{n-1}\cdot [1+(n+1)\cdot x].
\]
Hence, taking into account the independence polynomial of $K_{1}$, we see
that for any positive integer $k$, there is a well-covered graph $G$, namely
$G\in \{K_{1},K_{n}^{*},n\geq 1\}$, such that $I(G;x)$ has $-1/k$ as a root
and, in addition, all its roots are real.

\begin{figure}[h]
\setlength{\unitlength}{1cm}%
\begin{picture}(5,1.5)\thicklines

  \multiput(5,0)(1,0){3}{\circle*{0.29}}
  \multiput(5,1)(1,0){3}{\circle*{0.29}}
  \put(9,0){\circle*{0.29}}
  \put(9,1){\circle*{0.29}}
  \put(5,0){\line(0,1){1}}
  \put(6,0){\line(0,1){1}}
  \put(7,0){\line(0,1){1}}
  \put(5,0){\line(1,0){1}}
  \put(6,0){\line(1,0){1}}
  \put(9,0){\line(0,1){1}}
  \multiput(7,0)(0.125,0){16}{\circle*{0.07}}
  \put(4.65,0){\makebox(0,0){$b_{1}$}}
  \put(4.65,1){\makebox(0,0){$a_{1}$}}
  \put(6.3,0.3){\makebox(0,0){$b_{2}$}}
  \put(6.35,1){\makebox(0,0){$a_{2}$}}
  \put(7.3,0.3){\makebox(0,0){$b_{3}$}}
  \put(7.35,1){\makebox(0,0){$a_{3}$}}
  \put(9.37,0){\makebox(0,0){$b_{n}$}}
  \put(9.37,1){\makebox(0,0){$a_{n}$}}

 \end{picture}
\caption{The centipede $W_{n}$.}
\label{fig1}
\end{figure}
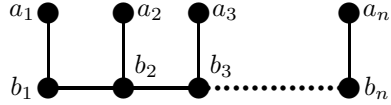

Let us consider the tree $W_{n}=P_{n}^{*},n\geq 1$, that we call a
\textit{centipede} (see\textit{\ }Figure \ref{fig1}).

In \cite{LevitMan2} it is noticed that for any $n\geq 2,I(W_{n};x)$
satisfies the recursion
\[
I(W_{n};x)=(1+x)\cdot (I(W_{n-1};x)+x\cdot
I(W_{n-2};x)),I(W_{0};x)=1,I(W_{1};x)=1+2x.
\]
In \cite{Arocha}, Arocha shows that $I(P_{n};x)=F_{n+1}(x)$, where
$F_{n}(x)$ are the so-called \textit{Fibonacci polynomials}, i.e.,
these polynomials are defined recursively by the following
formulae: $F_{0}(x)=1,F_{1}(x)=1,F_{n}(x)=F_{n-1}(x)+xF_{n-2}(x)$.
Based on this recurrence, one can deduce that
\[
I(P_{n};x)=\sum\limits_{j=0}^{\lfloor (n+1)/2\rfloor }{n+1-j
\choose j}\cdot x^{j}.
\]
Now, the equality $W_{n}=P_{n}^{*}$ and Theorem
\ref{th2}\emph{(i)} provide us with an
explicit form for the coefficients of $I(W_{n};x)=I(P_{n}^{*};x)=\sum%
\limits_{k=0}^{n}t_{k}\cdot x^{k}$, namely,
\[
t_{k}=\sum\limits_{j=0}^{k}{n-j \choose n-k}\cdot {n+1-j \choose
j},k\in \{0,1,2,...,n\}.
\]

The following result is a strengthening of Proposition \ref{prop1}\emph{(iv)}%
, as it concerns the real roots of the independence polynomial of a
well-covered graph.

\begin{proposition}
Let $G$ be a connected well-covered graph of girth $\geq 6$, which is not
isomorphic to $C_{7},K_{1},K_{2}$. Then the real roots of its independence
polynomial are in $[-1,-1/n)$, where $n=2\alpha (G)$.
\end{proposition}

\setlength {\parindent}{0.0cm}\textbf{Proof.} According to Theorem \ref{th3}%
, $G$ has a perfect matching consisting of only pendant edges, i.e., $%
G=H^{*} $ for some graph $H$. Hence, by Theorem \ref{th4}\emph{(vi)}, $%
I(G;x) $ has no real root $<-1$.\setlength
{\parindent}{3.45ex}

Further, Proposition \ref{prop1}\emph{(iv)} implies that any real root $%
x_{0} $ of $I(G;x)$ satisfies $\left| x_{0}\right| \geq 1/n$, while $1/n$ is
achieved only for a complete graph, i.e., only for $K_{2}$, in our case.
\rule{2mm}{2mm}\newline

Let us remark that\emph{\ }$I(K_{n};x),n\in \{1,2\}$, has a root at $-1/n$,
while not all the roots of $I(C_{7};x)=1+7x+14x^{2}+7x^{3}$ belong to $%
[-1,-1/7)$. More precisely, $I(C_{7};x)$ has at least one root in the
interval $(-2,-1)$, because $I(C_{7};-1)\cdot I(C_{7};-2)=-13$.

\begin{proposition}
\label{prop3}For any graph $G$ on $n\geq 2$ vertices, the following
assertions are valid:

\emph{(i)} $\max \{-\frac{\alpha \left( G\right) }{n},-\frac{1}{\omega (G)}%
\}\leq \xi _{\max }<-\frac{1}{2n-1}$;

\emph{(ii)} any complex root $z_{0}$ of $I(G;x)$ satisfies $\frac{1}{2n-1}%
<\left| z_{0}\right| $.
\end{proposition}

\setlength {\parindent}{0.0cm}\textbf{Proof.} As we saw in the
proof of Theorem \ref{th4}\emph{(iii)}, there is a bijection
\[
\begin{tabular}{c}
$\medskip f:A\rightarrow B,\ f(x)=\frac{x}{1-x},\ where$ \\
$A=\{x:I(G;x)=0\},B=\{x:I(G^{*};x)=0,x\neq -1\}.$%
\end{tabular}
\]
Now, if $x_{0}\in A$, then, according to Theorem \ref{th3}, the
corresponding root $f(x_{0})=\frac{x_{_{0}}}{1-x_{0}}\in B$ satisfies $%
\left| f(x_{0})\right| \geq 1/2n$.\setlength
{\parindent}{3.45ex} The equality $\left| f(x_{0})\right| =1/2n$ appears if $%
G^{*}$ is a complete graph, i.e., only for $G=K_{1}$. Taking now $G\neq
K_{1} $, we deduce that $\left| f(x_{0})\right| >1/2n$.

\emph{Case 1.} The root $x_{0}$ is real. In fact, $x_{0}<0$, and the
relation $\left| f(x_{0})\right| >1/2n$ leads to $-2nx_{0}>1-x_{0}$, which
gives $x_{0}<-1/(2n-1)$.

The relation $\max \{-\frac{\alpha \left( G\right) }{n},-\frac{1}{\omega (G)}%
\}\leq \xi _{\max }$ follows from Proposition \ref{prop1}\emph{(i),(v)}.

\emph{Case 2.} The root $z_{0}$ is not real. Then, $\left|
z_{_{0}}/(1-z_{0})\right| >1/2n$ implies
\[
2n\left| z_{_{0}}\right| >\left| 1-z_{0}\right| \geq \left| 1-\left|
z_{0}\right| \right| .
\]
Hence, $-2n\left| z_{_{0}}\right| <1-\left| z_{_{0}}\right| <2n\left|
z_{_{0}}\right| $ and further, $\left| z_{_{0}}\right| >1/(2n+1)$.

It is pretty amusing that one can not improve this bound using only simple
algebraic transformations. The proof of the bound $\frac{1}{2n-1}$ makes use
of Proposition\emph{\ }\ref{prop3}\emph{(i)} and Proposition \ref{prop1}%
\emph{(iii) }claiming that\emph{\ }a root of smallest modulus of $I(G;x)$ is
real, for any graph $G$. \rule{2mm}{2mm}

\begin{corollary}
If $T$ is a well-covered tree on $n\geq 4$ vertices,

then $-1=\xi _{\min }$ and $-\frac{1}{2}\leq \xi _{\max }<-\frac{1}{n}$.
\end{corollary}

\section{An application}

\begin{figure}[h]
\setlength{\unitlength}{1cm}%
\begin{picture}(5,1.2)\thicklines

  \multiput(7.2,0)(1,0){3}{\circle*{0.29}}
  \multiput(7.2,1)(1,0){3}{\circle*{0.29}}
  \put(7.2,0){\line(1,0){2}}
  \put(7.2,0){\line(0,1){1}}
  \put(7.2,0){\line(1,1){1}}
  \put(8.2,0){\line(1,1){1}}
  \put(6.7,0.5){\makebox(0,0){$G_{3}$}}

  \multiput(10.6,0)(1,0){2}{\circle*{0.29}}
  \multiput(10.6,1)(1,0){2}{\circle*{0.29}}
  \multiput(12.6,0)(0,1){2}{\circle*{0.29}}
  \put(10.6,0){\line(1,0){1}}
  \put(10.6,0){\line(0,1){1}}
  \put(10.6,0){\line(1,1){1}}
  \put(10.6,1){\line(1,0){1}}
  \put(11.6,0){\line(0,1){1}}
  \put(10.1,0.5){\makebox(0,0){$G_{4}$}}

  \multiput(3.9,0)(1,0){3}{\circle*{0.29}}
  \multiput(4.9,1)(1,0){2}{\circle*{0.29}}
  \put(3.9,0){\line(1,0){2}}
  \put(3.9,0){\line(1,1){1}}
  \put(4.9,0){\line(0,1){1}}
  \put(5.9,0){\line(0,1){1}}
  \put(3.6,0.5){\makebox(0,0){$G_{2}$}}

  \multiput(0.8,0)(1,0){3}{\circle*{0.29}}
  \multiput(1.8,1)(1,0){2}{\circle*{0.29}}
  \put(0.8,0){\line(1,0){2}}
  \put(0.8,0){\line(1,1){1}}
  \put(1.8,1){\line(1,0){1}}
  \put(2.8,0){\line(0,1){1}}
  \put(0.5,0.5){\makebox(0,0){$G_{1}$}}

  \end{picture}
\caption{Non-isomorphic ($G_{1},G_{2}$ are also well-covered) graphs having
the same independence polynomial $I(G_{1};x)=I(G_{2};x)$ and $%
I(G_{3};x)=I(G_{4};x)$.}
\label{fig63}
\end{figure}
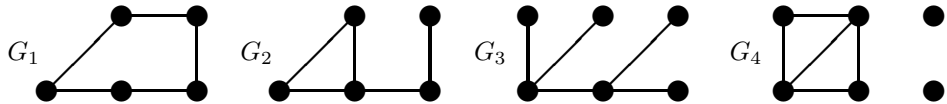

Let us observe that if $G$ and $H$ are isomorphic, then $I(G;x)=I(H;x)$. The
converse is not generally true. For instance, the graphs $%
G_{1},G_{2},G_{3},G_{4}$ depicted in Figure \ref{fig63} are non-isomorphic,
while $I(G_{1};x)=I(G_{2};x)=1+5x+5x^{2}$, and $%
I(G_{3};x)=I(G_{4};x)=1+6x+10x^{2}+6x^{3}+x^{4}$.

\begin{corollary}
\label{cor1}The following statements are true:

\emph{(i)} the graphs $G$ and $H$ are isomorphic if and only if $G^{*}$ and $%
H^{*}$ are isomorphic.

\emph{(ii)} $I(G;x)=I(H;x)$ if and only if $I(G^{*};x)=I(H^{*};x)$.
\end{corollary}

\setlength {\parindent}{0.0cm}\textbf{Proof.} \emph{(i)} The assertion
follows from the definition of $G^{*}$ and $H^{*}$, because any isomorphism $%
f:G\rightarrow H$ can be extended to an isomorphism $f^{*}:G^{*}\rightarrow
H^{*}$, while $f$ can be obtained as the restriction of an isomorphism $%
f^{*}:G^{*}\rightarrow H^{*}$ to $G$.\setlength
{\parindent}{3.45ex}

\emph{(ii)} Let $I(G;x)=I(H;x)$. Then $\alpha (G)=\alpha (H),\left|
V(G)\right| =\left| V(H)\right| =n$ and $I(G;x)=I(H;x)=\sum\limits_{k=0}^{%
\alpha (G)}s_{k}x^{k}$. According to Theorem \ref{th4}\emph{(i)}, it follows
that $I(G^{*};x)=\sum\limits_{k=0}^{\alpha (G)}s_{k}\cdot x^{k}\cdot
(1+x)^{n-k}=I(H^{*};x)$.

Conversely, assume that $I(G^{*};x)=I(H^{*};x)$. Hence, $\alpha
(G^{*})=\alpha (H^{*})=n$ and $\left| V(G)\right| =\left| V(H)\right| =n$.
According to Theorem \ref{th4}\emph{(ii)}, we infer that $%
I(G^{*};x)=(1+x)^{n}\cdot I(G;x/(1+x))$, and $I(H^{*};x)=(1+x)^{n}\cdot
I(H;x/(1+x))$. Therefore, the relation $I(G^{*};x)=I(H^{*};x)$ implies $%
I(G;x)=I(H;x)$. \rule{2mm}{2mm}\newline

Stevanovic \cite{Stevan} proved that the threshold graphs are clique-unique
graphs. It follows that the complements of threshold graphs are also
clique-unique graphs, since the class of threshold graphs is closed under
complement. Moreover, taking into account Corollary \ref{cor1}, we infer
that all the graphs of the family $\{G^{*}:G$ \textit{is a threshold graph}$%
\}$ are clique-unique graphs.

\begin{figure}[h]
\setlength{\unitlength}{1cm}%
\begin{picture}(5,1.2)\thicklines

  \multiput(1,0)(1,0){6}{\circle*{0.29}}
  \multiput(1,1)(2,0){2}{\circle*{0.29}}
  \multiput(4,1)(2,0){2}{\circle*{0.29}}
  \put(1,0){\line(1,0){5}}
  \put(1,0){\line(0,1){1}}
  \put(3,0){\line(0,1){1}}
  \put(4,0){\line(0,1){1}}
  \put(6,0){\line(0,1){1}}

  \put(0.5,0.3){\makebox(0,0){$T_{1}$}}

  \multiput(7.6,0)(1,0){6}{\circle*{0.29}}
  \multiput(9.6,1)(1,0){3}{\circle*{0.29}}
  \put(7.6,1){\circle*{0.29}}
  \put(7.6,0){\line(1,0){5}}
  \put(7.6,0){\line(0,1){1}}
  \put(9.6,1){\line(1,0){1}}
  \put(10.6,0){\line(0,1){1}}
  \put(11.6,0){\line(0,1){1}}

  \put(7.1,0.3){\makebox(0,0){$T_{2}$}}

 \end{picture}
\caption{Non-isomorphic trees with the same independence
polynomial.} \label{fig505}
\end{figure}
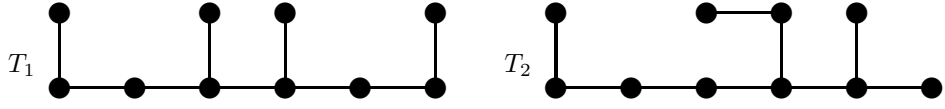

Recently, Dohmen, P\"{o}nitz and Tittmann \cite{DohPonTit} have found two
non-isomorphic trees having the same independence polynomial. These trees, $%
T_{1}$ and $T_{2}$, are depicted in Figure \ref{fig505}. They are clearly
non-isomorphic, while
\[
I(T_{1};x)=I(T_{2};x)=1+10x+36x^{2}+58x^{3}+42x^{4}+12x^{5}+x^{6}.
\]

Hence, according to Corollary \ref{cor1}, $I(T_{1}^{*};x)=I(T_{2}^{*};x)$,
while $T_{1}^{*},T_{2}^{*}$ are not isomorphic, because $T_{1},T_{2}$ are
not isomorphic. Moreover, $I((T_{1}^{*})^{*};x)=I((T_{2}^{*})^{*};x)$, while
$(T_{1}^{*})^{*},(T_{2}^{*})^{*}$ are not isomorphic. In this way, for any $%
k\geq 1$, we can find two non-isomorphic well-covered trees of size $10\cdot
2^{k}$, having the same independence polynomial. In addition, using the same
trees from Figure \ref{fig505}, it is easy to see that $T_{1}\cup
T_{1},T_{2}\cup T_{2}$ are not isomorphic, while $I(T_{1}\cup
T_{1};x)=I(T_{2}\cup T_{2};x)=I(T_{1};x)\cdot I(T_{2};x)$. Similarly, $%
T_{1}\cup T_{1}\cup T_{2},T_{1}\cup T_{2}\cup T_{2}$ are not isomorphic,
while $I(T_{1}\cup T_{1}\cup T_{2};x)=I(T_{1}\cup T_{2}\cup
T_{2};x)=I(T_{1};x)\cdot I(T_{2};x)\cdot I(T_{1};x)$ etc. Consequently, for
any $k\geq 1$, we can find two non-isomorphic well-covered forests of size $%
10\cdot k$, having the same independence polynomial.

In other words, the independence polynomial does not distinguish between
non-isomorphic trees. However, the following theorem claims that spiders are
uniquely defined by their independence polynomials in the context of
well-covered trees.

\begin{theorem}
\label{th8}The following statements are true:

\emph{(i)} for any $n\geq 2$, the independence polynomial of the
well-covered spider $S_{n}$ is
\[
I(S_{n};x)=(1+x)\cdot \left\{ 1+\sum\limits_{k=1}^{n}\left[ {n \choose k}%
\cdot 2^{k}+{n-1 \choose k-1}\right] \cdot x^{k}\right\} ;
\]

\emph{(ii) }if $G^{*}$ is connected, then the multiplicity of $-1$ as a root
of $I(G^{*};x)$ equals $1$

\qquad if and only if $G$ is isomorphic to $K_{1,n},n\geq 1$;

\emph{(iii)} if $G^{*}$ is connected, $I(G^{*};x)=I(T;x)$ and $T$ is a
well-covered spider,

\qquad then $G^{*}$ is isomorphic to $T$.
\end{theorem}

\setlength {\parindent}{0.0cm}\textbf{Proof.} \emph{(i)} If $G=K_{1,n},n\geq
2$, then $I(G;x)=\sum\limits_{k=0}^{\alpha (G)}s_{k}\cdot x^{k}=1+(n+1)\cdot
x+\sum\limits_{k=2}^{n}{n \choose k}\cdot x^{k}$ and $G^{*}=S_{n}$.%
\setlength
{\parindent}{3.45ex} Therefore, according to Theorem \ref{th2}, we obtain:
\begin{eqnarray*}
I(S_{n};x) &=&\sum\limits_{k=0}^{\alpha (G)}s_{k}\cdot x^{k}\cdot
(1+x)^{n+1-k}=\sum\limits_{k=0}^{n}s_{k}\cdot x^{k}\cdot (1+x)^{n+1-k} \\
&=&(1+x)^{n+1}+(n+1)\cdot x\cdot (1+x)^{n}+\sum\limits_{k=2}^{n}{n \choose k}%
\cdot x^{k}\cdot (1+x)^{n-k} \\
&=&(1+x)\cdot \left\{ x\cdot (1+x)^{n-1}+\sum\limits_{k=0}^{n}{n \choose k}%
\cdot x^{k}\cdot (1+x)^{n-k}\right\} \\
&=&(1+x)\cdot \left\{ 1+\sum\limits_{k=1}^{n}t_{k}\cdot x^{k}\right\} .
\end{eqnarray*}
Let us notice that the coefficient of $x^{k}$ is
\begin{eqnarray*}
t_{k} &=&{n-1 \choose k-1}+\sum\limits_{j=0}^{k}\left\{ {n \choose
j}\cdot
{n-j \choose k-j}\right\} \\
&=&{n-1 \choose k-1}+{n \choose k}\cdot \sum\limits_{j=0}^{k}{k \choose j}=%
{n-1 \choose k-1}+{n \choose k}\cdot 2^{k}.
\end{eqnarray*}

Consequently, $I(S_{n};x)=(1+x)\cdot \left\{
1+\sum\limits_{k=1}^{n}\left[ {n \choose k}\cdot 2^{k}+{n-1
\choose k-1}\right] \cdot x^{k}\right\} $, (for a different proof
of this relation, see \cite{LevitMan3}).

\emph{(ii) }According to Theorem \ref{th4}\emph{(vi)}, the multiplicity of $%
-1$ as a root of $I(G^{*};x)$ equals $\alpha (G^{*})-\alpha (G)=\left|
V(G)\right| -\alpha (G)$.

Now, if $-1$ is a simple root of $I(G^{*};x)$, then $\alpha (G)=\left|
V(G)\right| -1$, and because $K_{1,n}$ is the unique connected graph
satisfying this relation, it follows that $G$ is isomorphic to $K_{1,n}$.

Conversely, if $G$ is isomorphic to $K_{1,1}$, then $G^{*}=P_{4}$ and $%
I(P_{4};x)=1+4x+3x^{2}$ has $-1$ as a simple root.

Further, if $G$ is isomorphic to $K_{1,n},n\geq 2$, then, according to
Corollary \ref{cor1}, $G^{*}$ is isomorphic to $S_{n}$, and by Theorem \ref
{th4}\emph{(vi)},\emph{\ }$-1$ is a root of $I(G^{*};x)$ with the
multiplicity $\alpha (G^{*})-\alpha (G)=1$.

An alternative way to make the same conclusion is based on Theorem \ref{th8}%
\emph{(i)}. Since $I(G^{*};x)=I(S_{n};x)=(1+x)\cdot f(x)$, it follows that
\[
I(G^{*};1)=I(S_{n};1)=2\cdot \left\{ 1+\sum\limits_{k=1}^{n}\left[ {n \choose %
k}\cdot 2^{k}+{n-1 \choose k-1}\right] \right\} =2\cdot
(3^{n}+2^{n-1}).
\]
In other words, $f(1)=3^{n}+2^{n-1}$ is odd, and this ensures that $%
f(-1)\neq 0$, because, otherwise, if $f(-1)=0$, then $f(x)=(1+x)\cdot g(x)$,
and consequently, $f(1)=2\cdot g(1)$ is even. Therefore, $-1$ is a simple
root of $I(G^{*};x)$.

\emph{(iii) }Assume that $G^{*}$ is connected, $I(G^{*};x)=I(T;x)$ and $T$
is a well-covered spider.

If $T=K_{n},n=1,2$, then $I(G^{*};x)=1+nx$ and clearly $G^{*}$ is isomorphic
to $T$. If $T=P_{4}$, then $I(G^{*};x)=1+4x+3x^{2}$, and $G^{*}$ is
isomorphic to $P_{4}$, because there exists, by inspection, a unique
connected graph $H$ having $I(H;x)=1+4x+3x^{2}$, namely $P_{4}$. Further, if
$T=S_{n}=K_{1,n}^{*},n\geq 2$, then, the relation $I(G^{*};x)=I(T;x)$
implies, according to Theorem \ref{th8}\emph{(ii)}, that $I(G^{*};x)$ has $%
-1 $ as a simple root, and therefore, again by Theorem \ref{th8}\emph{(ii)},
$G^{*}$ is isomorphic to $T=S_{n}$. \rule{2mm}{2mm}\newline

Let us notice that the equality $I(G_{1};x)=I(G_{2};x)$ implies
\[
\left| V(G_{1})\right| =s_{1}=\left| V(G_{2})\right| \ and \
\left| E(G_{1})\right| =\frac{s_{1}^{2}-s_{1}}{2}-s_{2}=\left|
E(G_{2})\right| .
\]
Consequently, if $G_{1},G_{2}$ are connected, $I(G_{1};x)=I(G_{2};x)$ and
one of them is a tree, then the other must be a tree, as well. These
observations motivate the following conjecture.

\begin{conjecture}
\label{conj2}If $G$ is a connected graph and $T$ is a well-covered tree,
with the same independence polynomial, then $G$ is a well-covered tree.
\end{conjecture}

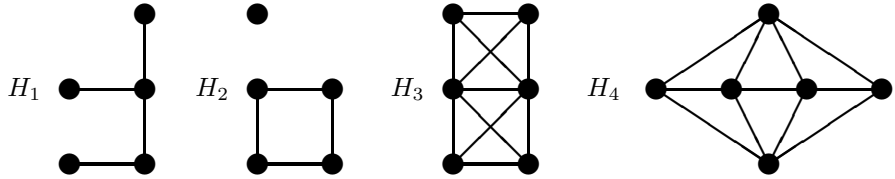
\begin{figure}[h]
\setlength{\unitlength}{1cm}%
\begin{picture}(5,2.3)\thicklines

  \multiput(1.2,0)(1,0){2}{\circle*{0.29}}
  \multiput(1.2,1)(1,0){2}{\circle*{0.29}}
  \put(2.2,2){\circle*{0.29}}
  \put(1.2,0){\line(1,0){1}}
  \put(1.2,1){\line(1,0){1}}
  \put(2.2,0){\line(0,1){2}}
  \put(0.6,1){\makebox(0,0){$H_{1}$}}

  \multiput(3.7,0)(0,1){3}{\circle*{0.29}}
  \multiput(4.7,0)(0,1){2}{\circle*{0.29}}
  \put(3.7,0){\line(1,0){1}}
  \put(3.7,0){\line(0,1){1}}
  \put(3.7,1){\line(1,0){1}}
  \put(4.7,0){\line(0,1){1}}
  \put(3.1,1){\makebox(0,0){$H_{2}$}}

 \multiput(6.3,0)(0,1){3}{\circle*{0.29}}
  \multiput(7.3,0)(0,1){3}{\circle*{0.29}}
  \put(6.3,0){\line(1,0){1}}
  \put(6.3,1){\line(1,0){1}}
  \put(6.3,2){\line(1,0){1}}
  \put(6.3,0){\line(0,1){2}}
  \put(7.3,0){\line(0,1){2}}
  \put(6.3,0){\line(1,1){1}}
  \put(6.3,1){\line(1,1){1}}
  \put(6.3,1){\line(1,-1){1}}
  \put(6.3,2){\line(1,-1){1}}
  \put(5.7,1){\makebox(0,0){$H_{3}$}}

  \put(10.5,0){\circle*{0.29}}
  \put(10.5,2){\circle*{0.29}}
  \multiput(9,1)(1,0){4}{\circle*{0.29}}
  \put(9,1){\line(1,0){3}}
  \put(9,1){\line(3,2){1.5}}
  \put(10,1){\line(1,2){0.5}}
  \put(9,1){\line(3,-2){1.5}}
  \put(10,1){\line(1,-2){0.5}}
   \put(12,1){\line(-3,2){1.5}}
   \put(11,1){\line(-1,2){0.5}}
  \put(12,1){\line(-3,-2){1.5}}
  \put(11,1){\line(-1,-2){0.5}}
  \put(8.3,1){\makebox(0,0){$H_{4}$}}

  \end{picture}
\caption{$I(H_{1};x)=I(H_{2};x)$ and $I(H_{3};x)=I(H_{4};x)$.}
\label{fig70}
\end{figure}

It is worth mentioning that changing the word ''tree'' for the word
''graph'' in Conjecture \ref{conj2} gives rise to a false assertion. For
example, $I(H_{1};x)=I(H_{2};x)=1+5x+6x^{2}+2x^{3}$, and $%
I(H_{3};x)=I(H_{4};x)=1+6x+4x^{2}$, where $H_{1},H_{2},H_{3},H_{4}$ are
depicted in Figure \ref{fig70}.

In other words, there exist a well-covered graph and a non-well-covered tree
with the same independence polynomial (e.g., $H_{2}$ and $H_{1}$), and also
a well-covered graph, different from a tree, namely $H_{4}$, satisfying $%
I(H_{3};x)=I(H_{4};x)$, where $H_{3}$ is not a well-covered graph.

\section{Conclusions}

To analyze the location structure of the roots of $I(G;x)$ in terms of
properties of $G$ seems to be a difficult task. For example, Hamidoune \cite
{Hamidoune} conjectures that for any \textit{claw-free} (i.e., a $K_{1,3}$%
-free) graph its independence polynomial has only real roots. Even in a
rather simple case of trees most of the relevant problems are open.

In this paper we found a number of properties concerning the interplay
between the (real) roots of the independence polynomials of graphs $G$ and $%
G^{*}$.

We also made an attempt to find some graph classes that can be defined by
their independence polynomials (\textit{independence-unique graphs}). In
this direction we succeeded in proving that well-covered spiders are
independence-unique among well-covered trees. If Conjecture \ref{conj2} is
true, than one may conclude that the independence polynomial distinguishes
between well-covered spiders and trees.

\end{document}